\documentclass[12pt]{amsart}

\newcommand {\Cee}    {{\mathbb  C}}

\newcommand {\fg}     {{\mathfrak{g}}}    %
\newcommand {\fgl}    {{\mathfrak{gl}}}  %
\newcommand {\fh}     {{\mathfrak{h}}}

\newcommand {\fp}    {{\mathfrak{p}}}
\newcommand {\fpo}    {{\mathfrak{po}}}
\newcommand {\fq}     {{\mathfrak{q}}}     

\newcommand {\fspe}   {{\mathfrak{spe}}}
\newcommand {\fpe}   {{\mathfrak{pe}}}
\newcommand {\fs}   {{\mathfrak{s}}}
\newcommand {\fsh}   {{\mathfrak{sh}}}

\newcommand {\fsq}    {{\mathfrak{sq}}}
\newcommand {\fsvect} {{\mathfrak{svect}}}

\newcommand {\fvect}  {{\mathfrak{vect}}}   %

\newcommand {\cal} {\mathcal}

\newcommand {\cO}     {{\cal O}}

\def \opname#1#2%
  {\expandafter\newcommand \csname #1\endcsname {{\mathop{#2}\nolimits}}}

\newcommand{\rmname}[1]
  {\expandafter\newcommand \csname #1\endcsname {{\operatorname{#1}}}}

\newcommand{\rmnameii}[2]
  {\expandafter\newcommand \csname #1\endcsname {{\operatorname{#2}}}}

\rmname{Ber}
\rmname{ber}
\rmname{Cliff}
\rmname{coind}
\rmname{diag}
\opname{Div}{div}
\rmname{End}
\rmname{id}
\rmname{ind}
\rmname{Ind}
\rmname{Inf}
\rmname{qet}
\rmname{qtr}
\rmname{ssym}
\rmname{str}
\rmname{tr}
\rmname{vpt}
\rmname{weyl}
\rmname{Weyl}
\rmname{Witt}
\opname{vvol}  {{v\hspace{-0.1ex}o\hspace{-0.02ex}l\/}}
\opname{Span} {{Span}}
\opname{Vol}  {{V\hspace{-0.55ex}o\hspace{-0.02ex}l\/}}

\rmname{Mat}
\rmname{Diff}
\rmname{Ker}
\rmname{Coker}
\rmname{Conn}
\rmname{Covect}
\rmname{Vect}
\rmname{Int}

\rmnameii {IM} {Im}
\rmnameii {RE} {Re}

\opname{Aut} {{A\hspace{-0.2ex}u\hspace{-0.1ex}t\/}}
\opname{GL} {{G\hspace{-0.3ex}L}}
\opname{SL} {{S\hspace{-0.3ex}L}}
\opname{Exp} {{E\hspace{-0.2ex}x\hspace{-0.1ex}p\/}}

\opname{U} {{U\/}}

\newcommand {\tto} {\longrightarrow}

\newcommand {\pderf}[2] {{\frac{\partial {#1}}{\partial {#2}}}}

\newcommand {\secno} {}
\newcommand {\ssecfont} {\normalfont\bf}

\newtheorem{Theorem}{\secno Theorem}

\newenvironment {th*}[1]
    {\gdef\thname{#1} \begin{thn}}%
    {\end{thn}}
\newtheorem{thn}[Theorem] {\thname}

\theoremstyle{definition}
\newtheorem{Example}{\secno Example}

\newenvironment {ex*}[1]
    {\gdef\thname{#1} \begin{exn}}%
    {\end{exn}}
\newtheorem{exn}[Theorem]{\thname}

\theoremstyle{remark}

\newtheorem{Remark}[Theorem]{\secno Remark}

\newenvironment {rem*}[1]
    {\gdef\thname{#1} \begin{remn}}%
    {\end{remn}}
\newtheorem{remn}[Theorem]{\thname}

\newcommand {\ssec}{\subsection*}

\newcommand {\ssbegin}[2]
  {\def \secno {\gdef \secno {}{\ssecfont #1. }}%
   \begin{#2}}
\begin{document}
	
\title{Casimir operators for Lie superalgebras}

\author{Dimitry Leites${}^1$, Alexander Sergeev${}^2$}

\address{${}^1$(Correspondence): Department of Mathematics, University
of Stockholm, Roslagsv.  101, Kr\"aftriket hus 6, SE-106 91,
Stockholm, Sweden; mleites@mathematik.su.se; ${}^2$Balakovo Inst.  of
Technique of Technology and Control, Branch of Saratov Technical
University, Balakovo, Saratov Region, Russia}

\thanks{The financial support of NFR is gratefully acknowledged.}

\begin{abstract}
Casimir operators --- the generators of the center of the enveloping
algebra --- are described for simple or close to them ``classical''
finite dimensional Lie superalgebras with nondegenerate symmetric even
bilinear form in: Sergeev A., The invariant polynomials on simple Lie
superalgebras.  Represent.  Theory 3 (1999), 250--280; math-RT/9810111
and for the ``queer'' series in: Sergeev A., The centre of enveloping
algebra for Lie superalgebra $Q(n,\,{C})$.  Lett.  Math.  Phys.  7,
no.  3, 1983, 177--179. 

Here we consider the remaining cases, and state conjectures proved for 
small values of parameter.

Under deformation (quantization) the Poisson Lie superalgebra
$\fpo(0|2n)$ on purely odd superspace turns into
$\fgl(2^{n-1}|2^{n-1})$ and, conjecturally, the lowest terms of the
Taylor series expansion with respect to the deformation parameter
(Planck's constant) of the Casimir operators for
$\fgl(2^{n-1}|2^{n-1})$ are the Casimir operators for $\fpo(0|2n)$. 

Similarly, quantization sends $\fpo(0|2n-1)$ into $\fq(2^{n-1})$ and 
the above procedure makes Casimir operators for $\fq(2^{n-1})$ into 
same for $\fpo(0|2n-1)$.

Casimir operators for the Lie superalgebra $\fvect(0|m)$ of vector
fields on purely odd superspace are only constants for $m>2$. 
Conjecturally, same is true for the Lie superalgebra $\fsvect(0|m)$ of
divergence free vector fields, and its deform, for $m>3$.

Invariant polynomials on $\fpo(0|2n-1)$ are also described.  They do
not correspond to Casimir operators.\end{abstract}



\maketitle 

This is a version of the paper published in: Ivanov E. et.  al. 
(eds.)  {\it Supersymmetries and Quantum Symmetries} (SQS'99, 27--31
July, 1999), Dubna, JINR, 2000, 409--411.  Meanwhile we have
understood how to prove our conjecture for $\fvect(0|m)$ and precisely
same idea was used by Noam Shomron \cite{Sh}, who additionally
described, to an extent, indecomposable $\fvect(0|m)$-modules. The 
only difference with the Dubna version is this and similar remarks.

As always, speaking about invariant polynomials on Lie superalgebras
(or their duals), one should bear in mind a more natural from super
point of view approach due to Shander \cite{Sha}: description of
nonpolynomial invariant functions.  Actually, even in non-super cases
people did consider non-polynomial elements, say, localizations of the
enveloping algebras, but there it does not lead to totally new
invariants.

\ssec{0.  Recapitulations} Recall that by the {\it Casimir elements}
for the Lie algebra or superalgebra $\fg$ we mean the generators of
the center $Z(\fg)$ of $U(\fg)$.  Their applications in description of
physical phenomena are well-known and numerous.

If $\fg$ possesses an invariant nondegenerate supersymmetric bilinear 
form $B$, then $\fg\simeq\fg^*$.  In this case, the description of the 
Casimir elements for $\fg$ is equivalent to the description of the 
algebra $I(\fg)$ of invariant polynomials on $\fg$ because $Z(\fg)$ is 
isomorphic to $I(\fg)$ (as vector spaces).  To describe $I(\fg)$ is, 
however, easier than $Z(\fg)$.  For the most complete description of 
$I(\fg)$ for simple finite dimensional Lie superalgebras and/or their 
``relatives'', i.e., close to them ``classical'' algebras, such as 
central extensions, etc., see \cite{S1}.  

This description of $I(\fg)$ does not cover $\fq(n)$, and $\fpo(0|m)$, 
and their relatives: it is a much more difficult problem to be 
solved by other means and ideas than those of \cite{S1}.

Thus, Casimir operators are described at the moment for 

$\bullet$ simple Lie superalgebras with Cartan matrix and their
``relatives'' \cite{S1}; (the first results are due to Berezin
\cite{B} and Kac \cite{K1, K2});

$\bullet$ $\fq(n)$ (but not for relatives), see \cite{S3}; 

$\bullet$ $\fpe(n)$ (no Casimir operators) and $\fspe(n)$, see
\cite{Sch} (implicitely) and \cite{S2} (for $n=3$); now we can refer
to an interesting paper of Maria Gorelik \cite{G}, where the
$\fspe(n)$-invariants are completely and explicitly described (though
the author herself does not like the result and did not even put on
arXive, to say nothing about publishing, we find it a very nice one;
it is based on her previous papers \cite{G2});

$\bullet$ $\fvect(0|m)$ (for $m>2$), see Shomron's paper \cite{Sh}.

There remain the following cases: 

$\bullet$ the divergence free sublagebras: the ``standard'' one, $\fsvect(0|m)$,
that preserves the volume element $\vvol$ and it deformation
$\widetilde{\fsvect(0|m)}$ preserving
$(1+\theta_{1}\cdot\theta_{m})\vvol$ for which,
conjecturally, there are no Casimir operators  (for $m>3$);

$\bullet$ $\fsq(n)$ and $\fp\fsq(n)$;

$\bullet$ $\fpo(0|m)$ considered below (but the general proof is still
lacking), $\fs\fpo(0|m)$ as well as $\fh(0|m)$ and $\fsh(0|m)$.

\begin{Example} For $\fg=\fgl(n)$ the Casimir operators 
corresponding to the first $n$ invariant polynomials $\tr(X^k)$, where 
$X\in \fg$, i.e., for $k=1$, \dots , $n$, generate the center of 
$U(\fg)$. The form $B$ is defined by the trace. \end{Example}

\begin{Example} The Lie superalgebra $\fgl(m|n)$ is a super analog of
$\fgl(n)$; the form $B$ is defined by the supertrace.  It is known
that the center of its universal enveloping algebra is NOT finitely
generated.  For the generators we may take $\str (X^k)$ for all
positive integer values of $k$.  (A comment: suppose we somehow proved
that apart from traces there are no invariant polynomilas on
$\fgl(n)$.  To prove that only finitely many of them generate the
algebra, consider the power series expansion of $\det(X-\lambda 1_n)$
with respect to $\lambda$.  The coefficients of $\lambda$ are
polynomials in $\tr(X^k)$; the characteristic polynomial yields a
recursive formula for $\tr(X^{n+1})$.  Contrarywise, the
superdeterminant, the Berezinian, is a rational function in
$\str(X^k)$, to express the Berezinian in terms of supertraces we need
all of them.)

On the other hand, any {\it rational} function on $\fgl(m|n)$ can be
expressed as a rational function of the first $n+m$ supertraces, cf. 
\cite{S1} with \cite{B} and \cite{K1, K2}).  \end{Example}

\begin{Example} The Lie superalgebra $\fpo(0|2n)$, the Poisson 
superalgebra, is one more super analog of $\fgl(n)$.  Indeed, there 
exists a well-known deformation (quantization) with parameter $\hbar$ 
such that at $\hbar=0$ we have $\fpo(0|2n)$ whereas at $\hbar\neq 0$ 
the Lie superalgebras from the parametric family are isomorphic to 
$\fgl(2^{n-1}|2^{n-1})$.

Recall, that the superspace of $\fpo(0|2n)$ is isomorphic to the 
Grassmann superalgebra on $2n$ generators $\xi$, $\eta$ and the 
bracket is given by the formula
$$
\{f, g\}=(-1)^{p(f)}\sum \left (\pderf{f}{\xi_{i}}\pderf{g}{\eta_{i}}+
\pderf{f}{\eta_{i}}\pderf{g}{\xi_{i}}\right ).
$$
The form $B$ on $\fpo(0|m)$ for any $m$ is given by the integral $B(f, 
g)=\int fg~\vvol$, where $\vvol$ is the volume element, and, clearly, 
the following are invariant polynomials on $\fpo(0|m)$:
$$
r_k=\int f^k\vvol(\xi, \eta),\text{ where } k=1, 2, \dots  \eqno{(*)}
$$
Observe, however, that, unlike 
$\fgl(m|n)$ case, these polynomials $r_k$ do not generate the whole 
algebra of invariant polynomials for $n\geq 2$.  A counterexample of 
least degree is given for $n=2$ by the Casimir operator whose radial part is 
equal to $x_1^2x_2^2(x_1^2-x_2^2)$.

Quantization turns the integral on functions that generate $\fpo(0|m)$
into the supertrace on supermatrices from $\fgl(2^{n-1}|2^{n-1})$ for
$m=2n$ and into queertrace on the general queer superalgebra
$\fq(2^{n-1})$ for $m=2n-1$.  This and other arguments, see \cite{BL},
indicates that $\fq(n)$ is one more super analog of $\fgl(n)$.  Recall
that $\qtr\left(\begin{pmatrix} A&B\cr B&A\end{pmatrix}\right)=\tr B$,
see \cite{BL}. \end{Example}

\begin{Remark} For any Lie superalgebra $\fg$ and its finite 
dimensional representation $\rho$, the polynomials
$$
P_{k, \rho}=\str \rho(X)^k
$$
are invariant ones.  For $\fgl(m|n)$ any degree $k$ homogeneous 
invariant polynomial can be represented as
$$
\mathop{\sum}\limits_{\rho} c_{\rho}
P_{k, \rho},       \eqno{(**)}
$$
where the sum runs over a finite number of finite dimensional 
representations, see \cite{S1}.  For Lie superalgebras $\fpo(0|2n)$ 
the situation is different: when $n\geq 2$ and $k\geq 6$ the 
polynomial $(*)$ can not be obtained as any finite sum of the form 
$(**)$, so we have strict inclusions $I(\fg)\supset 
I_{(*)}(\fg)\supset I_{(**)}(\fg)$, where $I_{(*)}$ and 
$I_{(**)}$ are algebras of polynomials of the form $(*)$ and $(**)$, 
respectively.  \end{Remark}

To describe the complete set of $\fpo(0|2n)$-invariants, recall that 
the Grassmann superalgebra $\Cee[\xi, \eta]$ can be deformed into the 
Clifford superalgebra $\Cliff_\hbar(\hat\xi, \hat \eta)$, where
$$
\hat\xi_{i}\hat \eta_{j}+\hat \eta_{j}\hat\xi_{i}=\delta_{ij}\hbar.
$$

For any subset $I=\{i_{1}< \dots < i_{l}\}\subset \{1, 2, \dots , n\}$ 
set $\xi_{I}=\xi_{i_{1}}\dots\xi_{i_{l}}$ and 
$\eta_{I}=\eta_{i_{1}}\dots\eta_{i_{l}}$; notations $\hat\xi_{I}$ and 
$\hat\eta_{I}$ are similar.  Define the linear map $Q: \Cee[\xi, 
\eta]\tto \Cliff_\hbar(\hat\xi, \hat \eta)$ by setting 
$Q(\xi_{I}\eta_{J})=\hat\xi_{I}\hat\eta_{J}$; before applying $Q$ to a 
monomial one has to reduce it to a normal form, say, to $\xi \eta$-form.

\ssbegin{4}{Lemma} $[Q(f), Q(g)]= \hbar Q(\{f, g\})+ \cO(\hbar^2)$,
where in the left hand side stands the supercommutator in the Clifford 
superalgebra, and $\{\cdot , \cdot\}$ in the right hand side is 
the Poisson bracket.
\end{Lemma} 
    
The Lie superalgebra associated with the associative superalgebra 
$\Cliff_\hbar(\hat\xi, \hat \eta)$ is $\fgl(2^{n-1}|2^{n-1})$; for it, 
the invariant polynomials are described above and in \cite{S1}.  The 
Lie superalgebra associated with the Clifford superalgebra generated 
by $2n-1$ elements is $\fq(2^{n-1})$ (see \cite{L}); for it, the 
Casimir elements are described in \cite{S3}.

\ssbegin{5}{Lemma} Let $F$ be an invariant polynomial on 
$\Cliff_\hbar$ considered as Lie superalgebra, i.e., on 
$\fgl(2^{n-1}|2^{n-1})$ or on its subalgebra $\fq(2^{n-1})$.  Let us 
represent it in the form $F=\mathop{\sum}\limits_{k\geq 
k_{0}}F_k\hbar^k$.  If $P\in I(\fpo(0|2n))$ is such that 
$Q(P)=F_{k_{0}}$, then $P$ is an invariant polynomial on $\fpo(0|m)$.
\end{Lemma} 

In other words, the lowest (with respect to $\hbar$-degree) components 
of invariant polynomials in the operators, i.e., on 
$\fgl(2^{n-1}|2^{n-1})$ or $\fq(2^{n-1})$, are invariant polynomials 
in their symbols --- the elements from $\fpo(0|m)$ for $m=2n$ or 
$2n-1$, respectively.  In particular, it is easy to see that
$$
\str(Q(f)^k) \quad (\text{or }\qtr(Q(f)^k))=\frac{1}{k}\left(\int 
f^k\vvol\right ) \hbar^n+\cO(\hbar^n).\eqno{(***)}
$$

\ssbegin{6}{Conjecture} {\em (\cite{S2})} All the invariant 
polynomials on $\fpo(0|m)$ for any $m$ (odd as well as even) may be 
obtained in the way indicated in Lemma 5.
\end{Conjecture}

This 15 years old conjecture is now proven for $m$ even and $\leq 
6$.  To prove it in full generality, we need new ideas as compared with 
those we know and successfully used in \cite{S1}: here the degree of 
complexity of computations performed \`a la \cite{S1} grows steeply 
with $m$.

\end{document}